\documentclass[12pt]{amsart}
\usepackage[top=1in, bottom=1in, left=1in, right=1in]{geometry}
\usepackage{dsfont}
\usepackage{amsmath} 
\usepackage{amsthm}
\usepackage{amssymb} 
\usepackage{mathtools}
\usepackage{times}
\usepackage[utf8]{inputenc}
\usepackage{indentfirst}
\usepackage{graphicx} 
\usepackage{color}
\usepackage[colorlinks=true]{hyperref}
\hypersetup{colorlinks=true, citecolor=green, linkcolor=green, filecolor=magenta, urlcolor=cyan}
\usepackage[shortlabels]{enumitem} 
\usepackage{mdframed}
\usepackage{tcolorbox}
\newtcolorbox{mybbox}{colback=blue!5!white,
colframe=blue!75!black}
\newtcolorbox{mygbox}{colback=green!5!white,
colframe=green!75!black}
\newtcolorbox{myrbox}{colback=red!5!white,
colframe=red!75!black}

\numberwithin{equation}{section}

\newtheorem{thm}{Theorem}[section]
\newtheorem{prop}[thm]{Proposition}
\newtheorem{lem}[thm]{Lemma}
\newtheorem{cor}[thm]{Corollary}
\theoremstyle{definition}

\theoremstyle{definition}

\theoremstyle{plain}
\newtheorem*{conjecture}{The Chowla Conjecture}
\theoremstyle{remark}

\newtheorem*{remark*}{Remark}

\newenvironment{feqn*}{\begin{mdframed}\begin{equation*}}{\vspace{1mm}
\end{equation*}\end{mdframed}}

\newcommand{\N}{\mathbb{N}}
\newcommand{\Z}{\mathbb{Z}}

\newcommand{\R}{\mathbb{R}}

\DeclareMathOperator{\supp}{supp}

\newcommand{\pthbgg}[1]{\bigg( #1 \bigg)}
\newcommand{\pthBg}[1]{\Big( #1 \Big)}

\newcommand{\abs}[1]{\lvert #1 \rvert}

\newcommand{\absbgg}[1]{\bigg\lvert #1 \bigg\rvert}
\newcommand{\absBg}[1]{\Big\lvert #1 \Big\rvert}

\newcommand{\bs}\boldsymbol{}
\newcommand{\psum}{\sideset{}{^{'}}\sum}
\renewcommand{\geq}{\geqslant}
\renewcommand{\leq}{\leqslant}

\newcommand{\eps}{\varepsilon}

\renewcommand{\mod}[1]{\,({\rm mod}\,#1)}

\definecolor{blue}{rgb}{.2,.6,.75}
\definecolor{green}{rgb}{.4,.7,.4}
\definecolor{red}{rgb}{1,0,0}

\begin{document}

\title[The Chowla conjecture and Landau-Siegel zeroes]{The Chowla conjecture and Landau-Siegel zeroes}

\author{Mikko Jaskari and Stelios Sachpazis}

\address{Department of Mathematics and Statistics, University of Turku, 20014 Turku, Finland}

\email{mikko.m.jaskari@utu.fi}

\address{Department of Mathematics and Statistics, University of Turku, 20014 Turku, Finland}
\email{stylianos.sachpazis@utu.fi}


\date{\today}

\begin{abstract}
Let $k\geq 2$ be an integer and let $\lambda$ be the Liouville function. Given $k$ non-negative distinct integers $h_1,\ldots,h_k$, the Chowla conjecture claims that $\sum_{n\leq x}\lambda(n+h_1)\cdots \lambda(n+h_k)=o(x)$ as $x\to\infty$. An unconditional answer to this conjecture is yet to be found, and in this paper, we take a conditional approach towards it. More precisely, we establish a non-trivial bound for the sums $\sum_{n\leq x}\lambda(n+h_1)\cdots \lambda(n+h_k)$ under the existence of a Landau-Siegel zero for $x$ in an interval that depends on the modulus of the character whose Dirichlet series corresponds to the Landau-Siegel zero. Our work constitutes an improvement over the previous related results of Germ\'{a}n and K\'{a}tai, Chinis, and Tao and Teräväinen.
\end{abstract}

\maketitle

\section{Introduction} 
The Liouville function $\lambda$ is the completely multiplicative function that is defined as $\lambda(n)=(-1)^{\Omega(n)}$ for all $n\in\N$, where $\Omega(n)=\sum_{p^{\alpha}\|n}\alpha$. The prime number theorem implies that
\begin{eqnarray*}\label{intro1}
\sum_{n\leq x} \lambda(n) = o(x)\quad \textrm {as}\quad x \to \infty,
\end{eqnarray*}
which means that the sign of $\lambda(n)$ changes frequently as $n$ grows. Chowla expected a more general version of this asymptotic to hold, and in 1965 \cite{chow}, he stated a conjecture which can be extended to the case of the Liouville function in the following way.
\begin{conjecture}
For distinct 
fixed integers $h_1, h_2,\dots, h_k\geq 0$, we have that
\begin{eqnarray*}\label{intro2}
\sum_{n \leq x}\lambda(n+h_1)\lambda(n+h_2)\cdots\lambda(n+h_k)=o(x)\quad \textrm {as}\quad x \to \infty.
\end{eqnarray*}
\end{conjecture}

Chowla's conjecture remains open till this day. So far, there has only been partial progress on it. Let us now take a look at the the main results that this progress evolved through. In 1985, Harman, Pintz and Wolke \cite{harman} studied the binary case with $h_1=0$ and $h_2=1$ and proved that
\begin{eqnarray*}\label{intro3}
-\frac{1}{3} + O\left(\frac{\log{x}}{x}\right) \leq \frac{1}{x}\sum_{n\leq x} \lambda(n)\lambda(n+1) \leq 1 - O_{\eps}\left(\frac{1}{(\log x)^{7+\eps}}\right)\quad(\eps>0).
\end{eqnarray*}
About 30 years later, in 2016, Matomäki and Radziwiłł \cite[Corollary 2]{short} improved this by showing that for any positive integer $h$, there exists $\delta(h) > 0$ such that
\begin{eqnarray*}\label{intro4}
\frac{1}{x}\left|\sum_{n \leq x}\lambda(n)\lambda(n+h)\right| \leq 1 - \delta(h)
\end{eqnarray*}
for all sufficiently large $x$. An averaged version of the conjecture was established by Matomäki, Radziwiłł and Tao \cite{averaged} in 2015. The three of them proved \cite[(1.3)]{averaged} that for every integer $k \geq 2$ and $H\in[10,x]$, we have
\begin{eqnarray*}\label{intro4}
\sum_{1\leq h_2,\dots,h_k\leq H}\left| \sum_{n \leq x} \lambda(n)\lambda(n+h_2)\cdots\lambda(n+h_k) \right| \ll k\left(\frac{\log\log{H}}{\log{H}}+\frac{1}{(\log x)^{1/3000}}\right)H^{k-1}x.
\end{eqnarray*}
Finally, a weaker version of Chowla's conjecture, the so-called logarithmically averaged Chowla conjecture which claims that
\begin{eqnarray*}\label{intro5}
\sum_{n \leq x} \frac{\lambda(n+h_1)\cdots\lambda(n+h_k)}{n} = o(\log{x})\quad \textrm{as}\quad x \to \infty,
\end{eqnarray*}
has been accomplished in stages for odd $k$ and $k=2$ in the works of Tao \cite{tao16}, Tao and Teräväinen \cite{tt18, tt19}, Helfgott and Radziwiłł \cite{helf}, and Pilatte \cite{pilatte}.

Despite the aforementioned partial advances towards Chowla's conjecture, the full conjecture still seems to be out of reach. To this end, one may wish to approach it conditionally, and in this paper, we are doing this by assuming the existence of Landau-Siegel zeroes. 

A \textit{Landau-Siegel zero}, often just called a \textit{Siegel zero}, is a real number $\beta$ associated to a real primitive Dirichlet character $\chi$ modulo $q$ such that $L(\beta,\chi)=0$ and 
$$\beta=1-\frac{1}{\eta\log q}$$
for some $\eta>0$. The number $\eta$ is called the \textit{quality} of the Landau-Siegel zero $\beta$ and from Siegel's theorem we have the ineffective bound
\begin{equation}\label{siegel}
    \eta\ll_{\eps}q^{\eps}
\end{equation}
for any $\eps>0$. It is widely believed that Siegel zeroes do not exist as their existence would come in contrast to the Generalized Riemann Hypothesis. However, the unlikely presence of Siegel zeroes has proven to have some interesting consequences. For example, a 1983 result of Heath-Brown \cite{heath} shows that if there exist infinitely many Siegel zeroes $\beta_j=1-1/(\eta_j\log q_j)$ such that $\eta_j\to +\infty$, then the Twin Prime Conjecture must be true.

The first result of the literature on Chowla's conjecture under Landau-Siegel zeroes comes from a work of Germán and Kátai from 2010. More precisely, in \cite[Theorem 2]{german}, they proved that there exists an absolute constant $c>0$ such that if $\chi$ is a primitive quadratic character modulo $q$ for which $L(\cdot,\chi)$ has a Landau-Siegel zero $\beta = 1 - 1/(\eta\log{q})$ with $\eta > \exp(\exp(30))$,
then
\begin{eqnarray*}
\frac{1}{x}\left|\sum_{n \leq x} \lambda(n)\lambda(n+1)\right| \leq \frac{c}{\log\log{\eta}}+o(1)
\end{eqnarray*}
for all $x\in[q^{10},q^{(\log\log{\eta})/3}]$. Under the same assumptions, Chinis \cite[Theorem 1.2]{jake} in 2021 extended the work of Germán and Kátai to the case of $k\geq 2$ distinct shifts and showed that
\begin{eqnarray*}
\frac{1}{x}\sum_{n\leq x} \lambda(n+h_1)\cdots\lambda(n+h_k)\ll_{k,h_1,\ldots,h_k} \frac{1}{(\log\log{\eta})^{1/2}(\log{\eta})^{1/12}}
\end{eqnarray*}
for all $x\in[q^{10},q^{(\log\log{\eta})/3}]$. A year later, in 2022, Tao and Teräväinen \cite[Corollary 1.8 (ii)]{tt} improved the bound of Chinis, as well as the $q$-range of $x$ that his bound holds for. In particular, under the existence of a Siegel zero $\beta=1-1/(\eta\log q)$ with $\eta\geq 10$, they showed that for a fixed $\eps\in(0,1)$, one has
\begin{eqnarray*}
\frac{1}{x}\sum_{n\leq x} \lambda(n+h_1)\cdots\lambda(n+h_k)\ll_{\eps,k,h_1,\ldots,h_k}\frac{1}{(\log{\eta})^{1/10}}
\end{eqnarray*}
for all $x$ in the interval $[q^{1/2+\eps},q^{\eta^{1/2}}]$. 

The purpose of this paper is to improve the result of Tao and Teräväinen in both directions. As we can see in Corollary \ref{corol} below, the following theorem, which is the main result of our work, leads to stronger bounds on the $k$-point correlations $\sum_{n\leq x} \lambda(n+h_1)\cdots\lambda(n+h_k)$ over a slightly wider $q$-range.
\begin{thm}\label{mainT}
Let $q\geq 2$ be a positive integer and let $\chi$ be a primitive quadratic character modulo $q$ such that $L(\cdot,\chi)$ has a real zero $\beta = 1 - 1/(\eta\log{q})$ with $\eta \geq 10$. We also fix an integer $k\geq 2$, distinct non-negative integers $h_1,\ldots,h_k$, and $\eps\in(0,1/2)$. There exists a constant $c=c(\eps,k)>0$ such that for $x=q^V$ with $V\in[1/2+\eps,\eta]$, we have
\begin{equation*}
\sum_{n\leq x} \lambda(n+h_1)\cdots\lambda(n+h_k)\ll\frac{xV}{\eta}+x\exp(-c\sqrt{V\log{\eta}}).
\end{equation*}
\end{thm}
\begin{cor}\label{corol}
Let $q, \eta, k$ and $\eps$ be as in the statement of Theorem \ref{mainT}. Let also $c$ be the constant from the bound of Theorem \ref{mainT} and consider distinct fixed integers $h_1,\ldots,h_k\geq 0$.  
\begin{enumerate}
    \item There exists a constant $c'=c'(\eps,k)>0$ such that
\begin{eqnarray*}
\sum_{n\leq x} \lambda(n+h_1)\cdots\lambda(n+h_k) \ll x\exp(-c'\sqrt{\log{\eta}})
\end{eqnarray*}
     for every $x \in [q^{1/2+\eps},q^{c^{-2}\log{\eta}}]$.
     \item For every $\delta\in(0,1)$, we have that
\begin{eqnarray*} 
\sum_{n\leq x} \lambda(n+h_1)\cdots\lambda(n+h_k) \ll_{\delta} \frac{x}{\eta^{1-\delta}}
\end{eqnarray*}
      for all $x \in (q^{c^{-2}\log{\eta}},q^{\eta^{\delta}}]$.
\end{enumerate}
\end{cor}
The implicit constants in Theorem \ref{mainT} and Corollary \ref{corol} depend on $\eps,k$ and the shifts $h_1,\ldots,h_k$, but we did not include them in the $\ll$ notation because the parameters $\eps,k,h_1,\ldots,h_k$ are all fixed in their statements. The same applies to several bounds occurring in the proof of Theorem \ref{mainT}. This way we alleviate the notation for the sake of a more convenient read.

We prove Theorem \ref{mainT} in Section \ref{sectionfive}. A major input in its proof is Lemma \ref{approx} which is a technical version of the fundamental lemma of sieve theory. Tao and Teräväinen did not use the beta-sieve which basically gives rise to Lemma \ref{approx}; Instead, they used a Selberg-type sieve which technically complicated the matters in their paper \cite{tt}. This difference between our work and theirs constitutes the main reason why we could improve their bound. Nonetheless, it could be possible to somehow alter their arguments and obtain bounds comparable to ours with a Selberg-type sieve.

The use of Lemma \ref{approx} as a main component in the improvement of the previous results was inspired by the work of Matomäki and Merikoski in \cite{mm}, where they assumed the existence of Siegel zeroes and refined a result of Tao and Teräväinen \cite{tt} on the binary ``Hardy-Littlewood"- type correlations $\sum_{n\leq x}\Lambda(n)\Lambda(n+h)$. As in our article, Matomäki and Merikoski, who were the ones who proved Lemma \ref{approx} in \cite{mm}, strengthened the result of Tao and Teräväinen by applying Lemma \ref{approx} in place of a Selberg-type sieve that the latter utilized for the more general ``Hardy-Littlewood-Chowla"-type correlations
$$\sum_{n\leq x}\Lambda(n+h_1)\cdots \Lambda(n+h_k)\lambda(n+h_1')\cdots\lambda(n+h_{\ell}')\quad (k\leq 2).$$ 

\subsection*{Notation} 
Let $r$ be a positive integer. For $a_1,\ldots,a_r\in\N$, we use the notation $(a_1,\ldots,a_r)$ for their greatest common divisor, and the notation $[a_1,\ldots,a_r]$ for their least common multiple.

Throughout the text, we denote the largest and smallest prime factor of an integer $n>1$ by $P^+(n)$ and $P^-(n)$, respectively. For $n=1$, we conventionally have $P^+(1)=0$ and $P^-(1)=+\infty$.

The capital Greek letters $\Phi$ and $\Psi$ also appear in the paper and they correspond to the classical counting functions of the rough and smooth integers, respectively. More precisely, for $x\geq y>1$, we have $\Phi(x,y)=\#\{n\leq x: P^-(n)>y\}$ and $\Psi(x,y)=\#\{n\leq x: P^+(n)\leq y\}$. 

Given two arithmetic functions $f$ and $g$, their Dirichlet convolution, denoted by $f\ast g$, is the arithmetic function which is defined as $(f\ast g)(n)=\sum_{ab=n}f(a)g(b)$ for all $n\in\N$. 

For $m\in\N$, the symbol $\tau_m$ denotes the $m$-fold divisor function given as $\tau_m(n)=\sum_{d_1\ldots d_m=n}1$ for all $n\in\N$. Finally, the lowercase Greek letters $\mu,\tau$ and $\varphi$ denote the M\"{o}bius, the divisor and the Euler functions, respectively.

\subsection*{Acknowledgements} We would like to thank Kaisa Matomäki for all her helpful comments and for drawing our attention to the problem addressed in this work. We also thank Cihan Sabuncu for the careful reading of an earlier version of the paper. During the making of this work, M. Jaskari was being supported by the Academy of Finland grant no. 346307 and the University of Turku Graduate School UTUGS. S. Sachpazis acknowledges support from the Academy of Finland grant no. 333707.

\section{The Core Ideas}
The purpose of this section is to explain why the existence of Siegel zeroes is a useful hypothesis when approaching Chowla's conjecture. Our explanation below also highlights the main steps in the proof of Theorem \ref{mainT}.  

If $L(\cdot,\chi)$ has a zero close to 1 for some Dirichlet character $\chi$ of modulus $q$, then $L(1,\chi)$ will be small by continuity. Consequently, assuming that Landau-Siegel zeroes exist, meaning that there exists a primitive quadratic character $\chi$ such that $L(\beta,\chi)=0$ for some real number $\beta$ close to 1, we expect that $L(1,\chi)^{-1}$ should be large. But,
$$L(1,\chi)^{-1}=\prod_p\pthbgg{1-\frac{\chi(p)}{p}},$$
so the Euler product of the right-hand side has to be large, too. In conjunction with Mertens' theorems, this suggests that $\chi(p)=-1=\lambda(p)$ for many large primes. This intuitive conclusion is the basis of arguments that use the presence of Landau-Siegel zeroes.

Now, let $z>1$ be some parameter. Following the work of Germán and Kátai \cite{german}, we define the completely multiplicative function $\lambda_z$ which is fully determined by its prime values
\begin{eqnarray}\label{lz}
\lambda_z(p)= 
\begin{cases} 
\lambda(p) &\text{if}\ p\leq z,\\ 
\chi(p) &\text{if}\ p>z.  
\end{cases}
\end{eqnarray}

According to the above discussion, in the presence of Siegel zeroes, the function $\lambda_z$ seems to be a good approximation to the Liouville function $\lambda$ for large $z$. Therefore, for a suitably chosen $z$, e.g. some small power of $x$, we expect
\begin{eqnarray}\label{N1}
\sum_{n\leq x} \lambda(n+h_1)\cdots\lambda(n+h_k)\approx \sum_{n\leq x} \lambda_z(n+h_1)\cdots\lambda_z(n+h_k).
\end{eqnarray}
In fact, such an approximation 
is provided by Lemma \ref{trans} which constitutes the first step in the proof of Theorem \ref{mainT}. 

Once the transition from the sums of $\lambda$ to the sums of $\lambda_z$ has been achieved with the introduction of an acceptable error, the study of the sums $\sum_{n\leq x}\lambda_z(n+h_1)\dots\lambda_z(n+h_k)$ is tractable. Let us elaborate on why is that. Without loss of generality, we assume that $h_1$ is the smallest shift. Then the right-hand side of (\ref{N1}) can be written as
\begin{eqnarray}\label{rew}
\sum_{\substack{a_1,\ldots,a_k\leq x,\\P^+(a_1\cdots a_k)\leq z}}\lambda(a_1)\cdots\lambda(a_k)\sum_{\substack{n_1\leq (x+h_1)/a_1,\\P^-(n_1\cdots n_k)>z,\\a_1n_1=a_in_i+h_1-h_i\\ \text{for all}\ i\in\{1,\ldots,k\}}}\chi(n_1)\cdots\chi(n_k).
\end{eqnarray}
The reason that makes the sums $\sum_{n\leq x}\lambda_z(n+h_1)\dots\lambda_z(n+h_k)$ easier to deal with lies on the innermost sum of (\ref{rew}), where we now have a primitive quadratic character $\chi$ instead of the Liouville function $\lambda$. This character $\chi$ possesses two crucial features that the Liouville function does not. The first one is its periodicity and the second one is the Weil bound (see Lemma \ref{Weil} below).

Let us now explain how the Weil bound is useful for the estimation of the expression in (\ref{rew}). First, we can approximate the indicator function $\mathds{1}_{P^-(\cdot)>z}$ by $1\ast w$ for some suitable sieve weight. This removes the condition $P^-(n_1\cdots n_k)>z$ from the innermost sum of (\ref{rew}) by introducing the terms $(1\ast w)(n_1\cdots n_k)$. After opening these convolutions, we can change the order of summation so that the innermost sums are character sums over polynomial values. So, provided that the contribution of the errors in the approximation $\mathds{1}_{P^-(\cdot)>z}\approx1\ast w$ can be handled, the estimation of (\ref{rew}) ends up being a simple application of the Weil bound to the new innermost character sums. The arguments of this explanation are rigorously executed in Section \ref{sectionfive}.

To summarize, the existence of Siegel zeroes is helpful because it puts our focus on the sums $\sum_{n\leq x} \lambda_z(n+h_1)\dots\lambda_z(n+h_k)$, which one can estimate by combining a sifting argument with the Weil bounds of $\chi$.

\section{Auxiliary Results}
In this section, we state and prove some preparatory results that are needed later for the proof of Theorem \ref{mainT}. We start with a lemma which will be used to bound the difference between the sums of the Liouville function $\lambda$ and those of the functions $\lambda_z$.
\begin{lem}\label{jo-ka}
Let $\chi$ be a primitive quadratic character modulo $q\geq 2$. Assume that $L(\cdot,\chi)$ has a real zero $\beta$ such that $\beta=1-1/(\eta\log q)$ for some $\eta\geq 10$. Let $z=q^v$ for some $v>0$. Then for any $x>z$, we have
\begin{eqnarray*}
\sum_{z<p\leq x}\frac{1+\chi(p)}{p}\ll
\begin{cases}
\displaystyle{\pthBg{\frac{1}{v^2\eta^{v/2}}+\frac{1}{z}}\pthBg{\frac{\log x}{\log z}}^3}& \text{for}\ v\in(0,2),\\
\displaystyle{\frac{\log x}{\eta\log q}}& \text{for}\ v\geq 2.
\end{cases}
\end{eqnarray*}
\end{lem}
\begin{proof}
The first branch of the estimate follows from Lemma 2.2 of \cite{mm}. Indeed, since $(1\ast\chi)(n)\geq 0$ for all $n\in\N$, for $v\in(0,2)$, Lemma 2.2 of \cite{mm} yields (Note that in \cite{mm}, the authors adopted the notation $\lambda=1\ast\chi$, whereas they denote the Liouville function by $\lambda_{\text{Liouville}}$.)
\begin{align*}
\sum_{z<p\leq x}\frac{1+\chi(p)}{p}&\leq \sum_{\substack{z<n\leq x\\P^-(n)>z}}\frac{(1\ast\chi)(n)}{n}\ll\pthBg{\frac{1}{v^2\eta^{v/2}}+\frac{v\log x}{\eta\log z}+\frac{1}{z}}\pthBg{\frac{\log x}{\log z}}^2\\
&\ll\pthBg{\frac{1}{v^2\eta^{v/2}}+\frac{v}{\eta}+\frac{1}{z}}\pthBg{\frac{\log x}{\log z}}^3\ll \pthBg{\frac{1}{v^2\eta^{v/2}}+\frac{1}{z}}\pthBg{\frac{\log x}{\log z}}^3.
\end{align*}

The second branch comes from Proposition 3.5 of \cite{tt}.
\end{proof}
The next lemma, which is found as Exercise 6(b) in \cite[end of Section III.5]{ten}, allows us to restrict our attention to the summands of $\sum_{n\leq x}\lambda_z(n+h_1)\cdots \lambda_z(n+h_k)$ for which the $z$-smooth parts $\prod_{p\leq z,\, p^{\alpha}\| n+h_i}p^{\alpha}$ of $n+h_i$ are relatively small for all $i\in\{1,\ldots,k\}$.
\begin{lem}\label{ste}
Consider the real numbers $u\geq 1$ and $x, z>1$. If $x\geq z^u$, then
$$\#\{n\leq x: \prod_{p\leq z,\, p^{\alpha}\| n}p^{\alpha}>z^u\}\ll xe^{-u/2}.$$
\end{lem}
\begin{proof}
We observe that  
\begin{eqnarray}\label{1}
\#\{n\leq x: \prod_{p\leq z,\, p^{\alpha}\| n}p^{\alpha}>z^u\}=\sum_{\substack{ab\leq x,\,z^u<a\\P^+(a)\leq z<P^-(b)}}1.
\end{eqnarray}

When $b=1$, the remaining sum over $a$ is bounded by $\Psi(x,z)$. When $b>1$, the condition $P^-(b)>z$ implies that $b>z$. So, for $b>1$, we have that $a\leq x/b<x/z$. For each such $a$ in (\ref{1}), there are at most $\Phi(x/a,z)$ integers $b$. Consequently, the relation (\ref{1}) gives
$$\#\{n\leq x: \prod_{p\leq z,\, p^{\alpha}\| n}p^{\alpha}>z^u\}\leq \Psi(x,z)+\sum_{\substack{z^u<a<x/z\\P^+(a)\leq z}}\Phi\pthBg{\frac{x}{a},z}.$$

For $t\geq z$, we have the bound $\Psi(t,z)\ll xe^{-\log t/(2\log z)}$ \cite[Theorem 1 in Section III.5.1]{ten}. Since $x/a>z$, we also have the standard estimate $\Phi(x/a,z)\ll x/(a\log z)$ which readily follows from Theorem 14.2 of \cite{dimb} with $f=\mathds{1}_{P^-(\cdot)>z}$. Therefore,
\begin{equation}\label{2}
\#\{n\leq x: \prod_{p\leq z,\, p^{\alpha}\| n}p^{\alpha}>z^u\}\ll xe^{-u/2}+\frac{x}{\log z}\sum_{\substack{z^u<a<x\\P^+(a)\leq z}}\frac{1}{a}.
\end{equation}

Using partial summation and the aforementioned bound on $\Psi(t,z)$ for $t\geq z$, it follows that the sum at the right-hand side of (\ref{2}) is $\ll \log z\cdot e^{-u/2}$. This concludes the proof of the lemma. 
\end{proof}
The following technical version of the fundamental lemma of sieve theory provides an approximation of the indicator function $\mathds{1}_{P^-(\cdot)>z}$ in terms of appropriate sieve weights. 
\begin{lem}\label{approx}
There exists a constant $\beta\geq 2$ such that the following holds for any $u\geq \beta$.

Let $z\geq 1$, and define 
$$z_r\vcentcolon=z^{((\beta-1)/\beta)^r}$$
\vspace{1mm}
for $r\in\N$. There exists an arithmetic function $w$ such that
\vspace{1mm}
\begin{itemize}
    \item $\abs{w(n)}\leq 1\quad \text{for all}\quad n\in\N$,
    \vspace{1mm}
    \item $\supp(w)\subseteq \{d\in\N:d\mid \prod_{p\leq z}p\text{ and } d\leq z^u\},\quad \text{and}$
    \item $\displaystyle{\mathds{1}_{P^-(n)>z}=(1\ast w)(n)+O\Big(\tau(n)^2\sum_{r\geq u-\beta}\mathds{1}_{P^-(n)>z_r}2^{-r}\Big)}\quad
    \text{for all}\quad n\in\N$.
\end{itemize}
\end{lem}
\begin{proof}
See the proof of \cite[Lemma 3.2(i)]{mm} with $A=1, \beta=\beta_1$ and a $u$ which here plays the role of the product $u\theta$ for some fixed $\theta\in(0,1/3)$.
\end{proof}
In the proof of Theorem \ref{mainT}, 
we also make use of the following extention of Shiu's theorem \cite[Theorem 1]{shiu} which is due to Henriot \cite{henr1, henr2}.
\begin{lem}\label{H}
Let $k$ be a positive integer and let $Q_1,\ldots,Q_k\in\Z[X]$ be $k$ pairwise coprime irreducible polynomials. Set $Q=Q_1\cdots Q_k$ and denote the degree and determinant of $Q$ by $g$ and $D$, respectively. For $n\in\N$ and $j\in\{1,\ldots,k\}$, we let $\rho_j(n)$ (respectively $\rho(n)$) denote the number of zeroes of $Q_j$ (respectively $Q$) modulo $n$. Assume that $Q$ has no fixed prime divisor and consider real numbers $A\geq 1$, $B\geq 1$, $\delta\in(0,1)$, and $\eps\in(0,1/(100g(g+1/\delta)))$. Suppose that $F:\N^k\to\R$ is a non-negative function such that
\begin{itemize}
\item $F(m_1n_1,\ldots,m_kn_k)=F(m_1,\ldots,m_k)F(n_1,\ldots,n_k)$\\
whenever $(m_1\cdots m_k,n_1\cdots n_k)=1$, and
\item $F(n_1,\ldots,n_k)\leq \min\{A^{\Omega(n_1\cdots n_k)},B(n_1\cdots n_k)^{\eps}\}\quad \text{for all}\quad (n_1,\ldots,n_k)\in\N^k.$
\end{itemize}
There exists a constant $c_0>0$, depending at most on $g, A, B$ and $\delta$, such that
\begin{eqnarray*}
\sum_{n\leq x}\!\!\!\!&F&\!\!\!\!\!(\abs{Q_1(n)},\ldots,\abs{Q_k(n)})\\
&\ll&\!\!\!x\prod_{p\mid D}\bigg(1+\sum_{\substack{0\leq \mu_j\leq \deg Q_j\\ \text{for all}\,j\in\{1,\ldots,k\},\\(\mu_1,\ldots,\mu_k)\neq (0,\ldots,0)}}F(p^{\mu_1},\ldots,p^{\mu_k})\bigg)\prod_{p\leq x}\bigg(1-\frac{\rho(p)}{p}\bigg)\\
&\times&\!\!\!\!\!\!\sum_{\substack{n_1\cdots n_k\leq x\\(n_1\cdots n_k,D)=1}}F(n_1,\ldots,n_k)\frac{\rho_1(n_1)\cdots \rho_k(n_k)}{n_1\cdots n_k}
\end{eqnarray*}
for all $x\geq c_0\Vert Q \Vert^{\delta}$, where $\Vert Q \Vert$ denotes the sum of the absolute values of the coefficients of $Q$. The implicit constant in the bound depends at most on $g, A, B$ and $\delta$.
\end{lem}
\begin{proof}
It follows from \cite[Theorem 3]{henr1} upon noticing that we can trivially bound the Euler product $\Delta_D$ in its statement by the simpler product
$$\prod_{p\mid D}\bigg(1+\sum_{\substack{0\leq \mu_j\leq \deg Q_j\\ \text{for all}\,j\in\{1,\ldots,k\},\\(\mu_1,\ldots,\mu_k)\neq (0,\ldots,0)}}F(p^{\mu_1},\ldots,p^{\mu_k})\bigg).\eqno\qedhere$$
\end{proof}
We will not need Lemma \ref{H} in its full generality. In fact, we can only think of $F$ as a product of divisor functions equipped with a coprimality condition, whereas $Q_1,\ldots,Q_k$ are simply going to be linear polynomials.

Another useful result for the proof of Theorem \ref{mainT} is the following generalization of the classical Weil bound for character sums \cite[Theorem 2C´]{sch}.
\begin{lem}\label{Weil}
Let $\chi$ be a primitive quadratic character modulo $q\geq 2$ and consider a polynomial $f\in\Z[X]$. Let also $M$ and $N\geq 0$ be two integers.  If $q^*$ is the product of the odd prime factors $p$ of $q$ for which $f$ is a constant multiple of a square polynomial modulo $p$, then
$$\sum_{M<n\leq M+N}\chi(f(n))\ll_{\eps}q^{\eps}\sqrt{q^*}\bigg(\frac{N}{\sqrt{q}}+\sqrt{q}\bigg).$$
\end{lem}
\begin{proof}
See the discussion preceding the statement of Lemma 3.7 in \cite[Subsection 3.4]{tt} 
\end{proof}
We close the section with an inequality which was stated without proof in \cite[Subsection 4.2]{jake}.
\begin{lem}\label{lcm}
Let $k$ be a positive integer. For given $m_1,\ldots,m_k\in\N$, we have
$$[m_1,\ldots,m_k]\geq\frac{m_1\cdots m_k}{\prod_{1\leq i<j\leq k}(m_i,m_j)}.$$
\end{lem}
\begin{proof}
For $k=1$ the lemma is trivial and we may assume that $k\geq 2$. 
We start by making use of the following elementary observation; For positive integers $a,b$ and $c$ for which $b\mid c$, it is true that $(a,b)\leq (a,c)$. We apply this property and obtain 
\begin{eqnarray}\label{1'}
(m_k,[m_1,\ldots,m_{k-1}])\leq (m_k,m_1\cdots m_{k-1}).
\end{eqnarray}

It is also easy to note that $(ab,c)\leq (a,c)(b,c)$ for $a,b,c\in\N$. 
Using this inequality, we obtain
\begin{equation}\label{2'}
(m_k,m_1\cdots m_{k-1})\leq (m_1,m_k)\cdots(m_{k-1},m_k).
\end{equation}

Now, we put (\ref{1'}) and (\ref{2'}) together and deduce that
$$[m_1,\ldots,m_k]=[m_k,[m_1,\ldots,m_{k-1}]]=\frac{m_k[m_1,\ldots,m_{k-1}]}{(m_k,[m_1,\ldots,m_{k-1}])}\geq \frac{m_k[m_1,\ldots,m_{k-1}]}{(m_1,m_k)\cdots(m_{k-1},m_k)}.$$

Based on this inequality, the proof is completed by induction on $k$.
\end{proof}
\section{Preparations for the Proof of Theorem \ref{mainT}} 
Throughout this and the following section, we assume that $\eps, k$ and $h_1, \dots, h_k$ are fixed. Hence, any dependence of the implicit constants on any of these numbers is ignored. Moreover, recall that for $z>1$, the completely multiplicative function $\lambda_{z}$ is defined by (\ref{lz}), where the Dirichlet character $\chi$ in (\ref{lz}) is the primitive real character modulo $q$ from the statement of Theorem \ref{mainT}. 
\subsection{Going from $\lambda$ to $\lambda_z$}
As was explained in Section 2, our first objective is to make a transition from the sums $\sum_{n \leq x} \lambda(n+h_1)\cdots\lambda(n+h_k)$ to the sums $\sum_{n \leq x} \lambda_{z}(n+h_1)\cdots\lambda_{z}(n+h_k)$. To this end, we prove the following.
\begin{lem}\label{trans}
Let $k,q\geq 2$ be natural numbers and let $h_1,\ldots,h_k$ be distinct non-negative integers. Let also $\chi$ be as in the statement of Theorem \ref{mainT} and write $z=q^v$ for some $v>0$. For $x>z$, we have that
\begin{eqnarray*}
\sum_{n \leq x}\bigg(\prod_{j=1}^k\lambda(n+h_j)-\prod_{j=1}^k\lambda_{z}(n+h_j)\bigg)\ll
\begin{cases}
\displaystyle{x\pthBg{\frac{1}{v^2\eta^{v/2}}+\frac{1}{z}}\pthBg{\frac{\log x}{\log z}}^3}& \text{for}\ v\in(0,2),\\
\displaystyle{\frac{x\log x}{\eta\log q}}&\text{for}\ v\geq 2.
\end{cases}
\end{eqnarray*}
\end{lem}
\begin{proof}
The beginning of the proof is similar to those of Lemma 3.1 in \cite{jake} and relation (3.3) in \cite[Section 3]{german}, but for the sake of completeness we write it down here. 

We will need the following inequality that one can prove by induction; If $m\in\N$ and $w_1,\ldots,w_m,$ $u_1,\ldots,u_m$ are complex numbers of modulus at most $1$, then 
\begin{equation}\label{auxineq}
\abs{w_1\cdots w_m-u_1\cdots u_m}\leq \sum_{i=1}^m\abs{w_i-u_i}. 
\end{equation}

Using this inequality, we infer that
\begin{align}\label{a}
\begin{split}
\sum_{n \leq x}\absBg{\prod_{j=1}^k\lambda(n+h_j) - \prod_{j=1}^k\lambda_{z}(n+h_j)}&\leq \sum_{n\leq x}\sum_{j=1}^k\abs{\lambda(n+h_j)-\lambda_z(n+h_j)}\\
&=\sum_{j=1}^k\sum_{n\leq x}\abs{\lambda(n+h_j)-\lambda_z(n+h_j)}\\
&=k\sum_{n \leq x} |\lambda(n)-\lambda_{z}(n)| + O(1).
\end{split}
\end{align}

In order to bound $\sum_{n \leq x} |\lambda(n)-\lambda_{z}(n)|$, we are exploiting the inequality (\ref{auxineq}) again. 
We combine it with the definition of $\lambda_z$ and get
\begin{align}\label{b}
\begin{aligned}
\sum_{n \leq x} |\lambda(n)-\lambda_{z}(n)| &\leq \sum_{n \leq x}\sum_{\substack{p^{\alpha}\|n \\ p>z}} |\lambda(p^{\alpha})-\chi(p^{\alpha})| \leq \sum_{\substack{p^{\alpha}\leq x \\ p>z}} |\lambda(p^{\alpha})-\chi(p^{\alpha})|\sum_{\substack{n \le x \\ p^{\alpha}|n}} 1 \\
&\ll x\sum_{\substack{p^{\alpha}\leq x \\ p>z}} \frac{|\lambda(p^{\alpha})-\chi(p^{\alpha})|}{p^{\alpha}}\leq x\sum_{z<p\leq x}\frac{1+\chi(p)}{p}+2x\sum_{p>z}\sum_{\alpha \geq 2}\frac{1}{p^{\alpha}} \\
&\ll x\sum_{z < p \leq x}\frac{1+\chi(p)}{p}+x\sum_{p>z}\frac{1}{p^2}\ll x\sum_{z < p \leq x}\frac{1+\chi(p)}{p}+\frac{x}{z}.
\end{aligned}
\end{align}

From (\ref{a}) and (\ref{b}), we deduce that
\begin{equation*}
\sum_{n \leq x}\bigg(\prod_{j=1}^k\lambda(n+h_j)-\prod_{j=1}^k\lambda_{z}(n+h_j)\bigg)\ll x\sum_{z < p \leq x}\frac{1+\chi(p)}{p}+\frac{x}{z}.
\end{equation*}

The proof is now completed by an application of Lemma \ref{jo-ka}. In the case $v<2$, Lemma \ref{jo-ka} finishes the proof immediately because $x/z$ is smaller than $x(\log x)^3/(z(\log z)^3)$ for $x>z$. When $v\geq 2$, the extra term $x/z$ is absorbed by the bound $x\log x/(\eta\log q)$ because in this case, we have the inequality $\log x>2\log q$ and the bound $\eta\ll z$ which follows from (\ref{siegel}).
\end{proof}
\subsection{Restricting the $z$-smooth parts}
Once Lemma \ref{trans} is used and the transition from the sums of $\lambda$ to the sums of $\lambda_z$ is completed, one's interest turns to the study of $\sum_{n \leq x} \lambda_{z}(n+h_1)\cdots\lambda_{z}(n+h_k)$. For technical reasons, we aim to exclude the terms of $\sum_{n \leq x} \lambda_{z}(n+h_1)\cdots\lambda_{z}(n+h_k)$ for which $\prod_{p\leq z,\, p^{\alpha}\| n+h_i}p^{\alpha}$ is relatively large for some $i\in\{1,\ldots,k\}$. To this end, in the current subsection, we assess the contribution stemming from these ``undesirable" summands.

Before doing that, we introduce the following notation for convenience; For $m\in\N$ and $z>1$, we write 
$$m_z \vcentcolon= \prod_{p^{\alpha} \| m, p \leq z} p^{\alpha}.$$

\begin{lem}\label{osz-sp}
Let $x\geq z>1$ and $u\geq 1$. If $x\geq z^u$, then
\begin{align*}
\sum_{n \leq x}\lambda_{z}(n+h_1)\cdots\lambda_z(n+h_k)=\sum_{\substack{n \leq x, \\ (n+h_i)_{z} \leq z^{u}\\ \text{for all } i \in \{1,\dots,k\}}} \lambda_{z}(n+h_1)\cdots\lambda_{z}(n+h_k)+O(xe^{-u/2}+1).
\end{align*}
\end{lem}
\begin{proof}
We readily see that 
\begin{align}\label{start}
\begin{split}
\sum_{n \leq x}\lambda_{z}(n+h_1)\cdots\lambda_z(n+h_k)&=\sum_{\substack{n \leq x, \\ (n+h_i)_{z} \leq z^{u}\\ \text{for all } i \in \{1,\dots,k\}}} \lambda_{z}(n+h_1)\cdots\lambda_{z}(n+h_k)\\
&+\sum_{\substack{n \leq x, \\ (n+h_i)_{z} > z^{u} \\ \text{for some } i \in \{1,\dots,k\}}} \lambda_{z}(n+h_1)\cdots\lambda_{z}(n+h_k).
\end{split}
\end{align}

Note that
\begin{align*}
\absbgg{\sum_{\substack{n \leq x, \\ (n+h_i)_{z} > z^{u} \\ \text{for some } i \in \{1,\dots,k\}}} \lambda_{z}(n+h_1)\cdots\lambda_{z}(n+h_k)}&\leq \sum_{i=1}^{k} \sum_{\substack{n \leq x \\ (n+h_i)_{z} > z^{u}}} 1 \leq \sum_{i=1}^{k} \sum_{\substack{h_i < m \leq x+h_i \\ m_{z} > z^{u}}} 1\\
&= \sum_{i=1}^{k} \#\{ m \leq x: \prod_{p^{\alpha}\| m, p\leq z} p^{\alpha} > z^u \} + O(1).
\end{align*}

Therefore, by using Lemma \ref{ste}, we deduce that 
\begin{equation}\label{*}
\sum_{\substack{n \leq x, \\ (n+h_i)_{z} > z^{u} \\ \text{for some } i \in \{1,\dots,k\}}} \lambda_{z}(n+h_1)\cdots\lambda_{z}(n+h_k) \ll xe^{-u/2} + 1.  
\end{equation}

The proof now follows by inserting (\ref{*}) into (\ref{start}).
\end{proof}
\subsection{Rewriting the ``main term" in Lemma \ref{osz-sp}} Let $x,z$ and $u$ be as in the statement of Lemma \ref{osz-sp}. We are going to rewrite the sum
\begin{equation}\label{defS1}
S\vcentcolon=\sum_{\substack{n \leq x, \\ (n+h_i)_{z} \leq z^{u}\\ \text{for all } i \in \{1,\dots,k\}}} \lambda_{z}(n+h_1)\cdots\lambda_{z}(n+h_k)   
\end{equation}
in a form that we use in the upcoming section. For the sequel, we assume without loss of generality that $h_1 = \min_{1 \leq i \leq k} h_i$. We write $n+h_i=a_in_i$ with $P^+(a_i)\leq z$ and $P^-(n_i)>z$, and so
\begin{align}\label{S1a}
S = \sum_{\substack{a_1, \dots, a_k \leq z^{u}, \\ P^{+}(a_1\cdots a_k) \leq z, \\ (a_i,a_j) |(h_i-h_j) \\ \text{for all } i,j \in \{1,\dots,k\} \textrm{ with } i\not=j}}\lambda(a_1\cdots a_k)\sum_{\substack{n_1 \leq \frac{x+h_1}{a_1}, \\ P^{-}(n_1\cdots n_k) > z, \\ a_{1}n_{1} = a_{i}n_{i}+h_1-h_i \\ \text{for all } i \in \{2,\dots,k\}}} \chi(n_1\cdots n_k).
\end{align}
We denote the inner sum in (\ref{S1a}) by $\Sigma$, that is,
\begin{align*}
\Sigma \vcentcolon= \sum_{\substack{n_1 \leq \frac{x+h_1}{a_1}, \\ P^{-}(n_1\cdots n_k) > z, \\ a_{1}n_{1} = a_{i}n_{i}+h_1-h_i \\ \text{for all } i \in \{2,\dots,k\}}} \chi(n_1\cdots n_k).
\end{align*}
Upon writing $h_i^* = h_i - h_1$ for $i \in \{1,\dots,k\}$, we see that
\begin{align}\label{D1}
\Sigma = \sum_{\substack{n \leq \frac{x+h_1}{a_1}, \\ P^{-}\left( \prod_{i=1}^{k}\frac{a_{1}n+h_{i}^*}{a_i} \right)>z, \\ a_{1}n\equiv -h_i^* \mod{a_i} \\ \text{for all } i \in \{2,\dots,k\}}} \chi\bigg( \prod_{i=1}^{k}\frac{a_{1}n+h_{i}^*}{a_i} \bigg).
\end{align}
Since we have the condition $(a_i,a_j)|(h_i-h_j)$ in the outermost sum of (\ref{S1a}), for every $i \in \{1,\dots,k\}$, one has that
\begin{align}\label{lincon}
\begin{split}
a_{1}n \equiv -h_i^* \mod{a_i} &\Longleftrightarrow \frac{a_1}{(a_1,a_i)}n \equiv -\frac{h_i^*}{(a_1,a_i)} \left(\mathrm{mod}\,\frac{a_i}{(a_1,a_i)}\right) \\
&\Longleftrightarrow n \equiv -\overline{a_{1,i}}\frac{h_i^*}{(a_1,a_i)}(\mathrm{mod}\, a_i^*),
\end{split}
\end{align}
where $a_i^* := a_{i}(a_1,a_i)^{-1}$, $a_{1,i}:=a_{1}(a_1,a_i)^{-1},$ and $\overline{a_{1,i}}$ is the inverse of the latter $(\mathrm{mod}\, a_i^*)$. Now, if the system of linear congruences (\ref{lincon}) is not soluble, then $\Sigma=0$ and there is no contribution from such systems. Consequently, we may limit the outermost sum of (\ref{S1a}) to those $a_1,\ldots,a_k$ for which the system of linear congruences (\ref{lincon}) is soluble. In this case, by the Chinese Remainder Theorem, there exists an integer $r^*$ such that
\vspace{1mm}
\begin{itemize}
\item $0 < r^* \leq [a_{2}^*, \dots, a_{k}^*],$
\item $n = \ell[a_{2}^*, \dots, a_{k}^*]+r^*-[a_{2}^*, \dots, a_{k}^*]$ for $\ell \in \N,$ and
\item $a_{1}r^* \equiv -h_i^* \mod{a_i}$ for all $i \in \{1,\dots,k\}$.
\end{itemize}
\vspace{1mm}
Applying these to (\ref{D1}), we derive that
\begin{align}\label{D2}
\Sigma = \sum_{\substack{l \leq \frac{x+h_1-a_{1}r^*}{a_{1}[a_{2}^*, \dots, a_{k}^*]}+1 \\ P^{-}(Q(\ell))>z}}\chi(Q(\ell)),
\end{align}
where
\begin{align*}
Q(X) := \prod_{i=1}^{k}\left( a_{1,i}\frac{[a_{2}^*, \dots, a_{k}^*]}{a_i^*}X+\frac{a_{1}r^*+h_{i}^*}{a_i}-a_{1,i}\frac{[a_{2}^*, \dots, a_{k}^*]}{a_{i}^*} \right)
\end{align*}
defines a polynomial of $\Z[X]$. 
Since $\Sigma$ denotes the inner sum of (\ref{S1a}), by inserting (\ref{D2}) into (\ref{S1a}), we obtain that
\begin{align}\label{S1b}
S = \psum_{\substack{a_1, \dots, a_k \leq z^{u}, \\ P^{+}(a_1\cdots a_k) \leq z, \\ (a_i,a_j) |(h_i-h_j) \\ \text{for all } i,j \in \{1,\dots,k\} \text{ with } i\not=j}}\lambda(a_1\cdots a_k)\sum_{\substack{l \leq \frac{x+h_1-a_{1}r^*}{a_{1}[a_{2}^*, \dots, a_{k}^*]}+1 \\ P^{-}(Q(\ell))>z}}\chi(Q(\ell)),
\end{align}
where the prime $^{'}$ indicates that the sum is taken over the integers $a_1,\ldots,a_k$ for which the system of congruences (\ref{lincon}) is soluble. 

Relation (\ref{S1b}) provides the form of $S$ that we will exploit.

\section{Proof of Theorem \ref{mainT}}\label{sectionfive}
In this section, we establish Theorem \ref{mainT}. To achieve this, we bound $S$ by starting from (\ref{S1b}), and then, once the estimation of $S$ is complete, we combine the definition (\ref{defS1}) of $S$ with Lemmas \ref{trans} and \ref{osz-sp}.

Let $z=q^v$ with
\begin{equation}\label{vch}
v=\min\Bigg\{\sqrt{\frac{V}{\log\eta}},\,2\Bigg\}.
\end{equation}

We now assume that $\eta$ is sufficiently large, otherwise $xV/\eta\gg x$, in which case, Theorem \ref{mainT} would immediately follow from the trivial bound $x$. Since $\eta$ is assumed to be large enough, Siegel's theorem, that is (\ref{siegel}), implies that $q$ is also sufficiently large and that $v\gg1/\sqrt{(\log q)}$. Therefore, $v \geq \log{k}/\log{q}$, and so $z = q^v \geq k$.



In (\ref{S1b}), we have that $P^+(a_1\cdots a_k)\leq z<P^-(Q(\ell))$, and so the inner sum $\Sigma$ of $S$ is taken over integers $\ell$ for which $(Q(\ell),a_1\cdots a_k)=1$. Moreover, since $z\geq k$, we also have that $P^-(Q(\ell))>z\geq k$, and this imposes the condition $(Q(\ell),\prod_{p\leq k}p)=1$ on $\Sigma$. According to the above, we can add the supplementary condition $(Q(\ell),a_1\cdots a_k\prod_{p \leq k} p)=1$ to the inner sum $\Sigma$ of (\ref{S1b}). Therefore,
\begin{align}\label{Sig}
\Sigma = \sum_{\substack{l \leq \frac{x+h_1-a_{1}r^*}{a_{1}[a_{2}^*, \dots, a_{k}^*]}+1 \\ P^{-}(Q(\ell))>z\\(Q(\ell),a_1\cdots a_k\prod_{p \leq k} p)=1}}\chi(Q(\ell)).
\end{align} 

We now detect the sifting condition $P^-(Q(\ell))>z$ in (\ref{Sig}) by applying Lemma \ref{approx} with $\beta=25k$ and $u=\eps V/(30kv)$.
Since $V \geq \frac{1}{2}+\varepsilon$ and $\eta$ is large, we observe that $u=\eps V/(30vk)\geq 50k$.
Hence,
\begin{align}\label{D3}
\begin{split}
\Sigma =&\, \sum_{\substack{d \leq z^{u}\\ d\mid\prod_{p\leq z}p \\ (d,a_1\cdots a_k\prod_{p \leq k} p)=1}}w(d)\sum_{\substack{\ell \leq \frac{x}{a_1[a_2^*,\dots,a_k^*]} \\ d|Q(\ell) \\ (Q(\ell),a_1\cdots a_k\prod_{p \leq k} p)=1}}\chi(Q(\ell)) \\
&+O\Bigg(z^u+ \sum_{r \geq \frac{u}{2}}2^{-r}\sum_{\substack{\ell \leq \frac{x}{a_1[a_2^*,\dots,a_k^*]} \\ P^{-}(Q(\ell))>z_r \\ (Q(\ell),a_1\cdots a_k\prod_{p \leq k} p)=1}}\tau(Q(\ell))^{2} \Bigg),
\end{split}
\end{align}
where $z_r = z^{((25k-1)/(25k))^r}$ and $w$ is the arithmetic function from the statement of Lemma \ref{approx}.
\subsection{Evaluation of the ``main term" of $\Sigma$ in (\ref{D3})}
We seek a bound for
\begin{align*}
\Sigma' := \sum_{\substack{d \leq z^{u}\\ d\mid\prod_{p\leq z}p \\ (d,a_1\cdots a_k\prod_{p \leq k} p)=1}}w(d)\sum_{\substack{\ell \leq \frac{x}{a_1[a_2^*,\dots,a_k^*]} \\ d|Q(\ell) \\ (Q(\ell),a_1\cdots a_k\prod_{p \leq k} p)=1}}\chi(Q(\ell)).
\end{align*}

We spot the condition $(Q(\ell),a_1\cdots a_k\prod_{p \leq k} p)=1$ using Möbius inversion, and so we obtain
\begin{align}\label{Dp1}
\begin{split}
\Sigma' &= \sum_{\substack{d \leq z^{u},\, d\mid\prod_{p\leq z}p \\ (d,a_1\cdots a_k\prod_{p \leq k} p)=1 \\ d'|a_1\cdots a_k\prod_{p \leq k} p}}w(d)\mu(d')\sum_{\substack{\nu \mod{dd'} \\ Q(\nu) \equiv 0 \mod{dd'}}}\sum_{\substack{\ell \leq \frac{x}{a_1[a_2^*,\dots,a_k^*]} \\ \ell \equiv \nu \mod{dd'} }} \chi(Q(\ell)) \\
&= \sum_{\substack{d \leq z^{u},\, d\mid\prod_{p\leq z}p \\ (d,a_1\cdots a_k\prod_{p \leq k} p)=1 \\ d'|a_1\cdots a_k\prod_{p \leq k} p}}w(d)\mu(d')\sum_{\substack{\nu \mod{dd'} \\ Q(\nu) \equiv 0 \mod{dd'}}}\sum_{m \leq \frac{x}{dd'a_1[a_2^*,\dots,a_k^*]}} \chi(R_{d,d',\nu}(m)),
\end{split}
\end{align}
where $R_{d,d',\nu}\in\Z[X]$ with $R_{d,d',\nu}(X):=Q(dd'X+\nu)$.

We want to apply Lemma \ref{Weil} to bound the innermost character sum in (\ref{Dp1}). In order to do this with $q^*\ll (dd',q)$, we have to verify that $R_{d,d',\nu}$ is not a constant multiple of a square polynomial modulo $p$ for those prime factors $p$ of $q$ for which $(p,dd')=1$ and $p > \max_{1 \leq i \leq k} h_i$. We do this by considering the following complementary cases; (i) $p \nmid a_i$ for all $i \in \{1,\dots,k\}$, and; (ii) $p | a_i$ for some $i \in \{1,\dots,k\}$.

We first assume that $p \nmid a_i$ for all $i \in \{1,\dots,k\}$. In this case, if $R_{d,d',\nu}$ is a square polynomial modulo $p$, its roots have even multiplicity, and this implies that there exist distinct indices $i,j \in \{1,\dots,k\}$ such that
\begin{align*}
&\,\overline{dd'}\cdot\overline{\,\,\frac{a_1[a_2^*,\dots,a_k^*]}{a_i}\,\,}\cdot\left( \frac{a_1r^*+h_i^*}{a_i} + \frac{(\nu-1)a_1[a_2^*,\dots,a_k^*]}{a_i} \right) \nonumber \\
\equiv&\,\overline{dd'}\cdot\overline{\,\,\frac{a_1[a_2^*,\dots,a_k^*]}{a_j}\,\,}\cdot\left( \frac{a_1r^*+h_j^*}{a_j} + \frac{(\nu-1)a_1[a_2^*,\dots,a_k^*]}{a_j} \right) \mod{p}.
\end{align*}
We multiply both sides of this congruence by $dd'(a_1[a_2^*,\dots,a_k^*])^{2}a_i^{-1}a_j^{-1}$ and infer that 
\begin{equation*}
\frac{a_1[a_2^*,\dots,a_k^*]}{a_ia_j}(h_i^* - h_j^*) \equiv 0 \mod{p}.
\end{equation*}
Since $p \nmid a_i$ for all $i\in\{1,\dots,k\}$, we get $h_i^* \equiv h_j^*\mod{p}$, which is equivalent to $h_i \equiv h_j \mod{p}$. But, since $h_i,h_j \leq \max_{1 \leq \kappa \leq k} h_\kappa < p$, we conclude that $h_i = h_j$, which is a contradiction for $i \not= j$ as the shifts $h_1,\ldots,h_k$ are distinct. Hence, $R_{d,d',\nu}$ is not a constant multiple of a square polynomial modulo $p$ in this case.

We now assume that $p | a_i$ for some $i \in \{1,\dots,k\}$. This implies that $p \nmid a_j$ for all $j \not= i$, otherwise $p|(a_i,a_j)$ for some $j\not=i$, which, combined with the condition $(a_i,a_j)|(h_i-h_j)$ of (\ref{S1b}), would imply the divisibility relation $p|(h_i-h_j)$ which contradicts the condition $p > \max_{1 \leq i \leq k} h_i$ that we are working under. By the same argument, if $i\neq 1$, then we can deduce that $p\mid a_i^*$. Combining both our conclusions, we observe that $p$ divides $a_{1,j}[a_2^*,\dots,a_k^*]{a_{j}^*}^{-1}=a_1[a_2^*,\dots,a_k^*]{a_j}^{-1}$ for $j\neq i$. Note also that 
$$a_{1,i}[a_2^*,\dots,a_k^*]{a_{i}^*}^{-1}\Big| \prod_{\substack{j=1, j \not= i}}^{k} a_j$$
because $[a_2^*,\dots,a_k^*]$ divides $a_2^*\cdots a_k^*$, which in turn divides $a_2\cdots a_k$. So, $p\nmid a_{1,i}[a_2^*,\dots,a_k^*]{a_{i}^*}^{-1}$, whereas $p\mid a_{1,j}[a_2^*,\dots,a_k^*]{a_{j}^*}^{-1}$ for $j\neq i$. We thus see that $R_{d,d',\nu}(X)\equiv c_1X+c_0 \mod{p}$ for some $c_0,c_1 \in \N$. This shows that $R_{d,d',\nu}$ is not a constant multiple of a square polynomial modulo $p$ in this case either.

Now that we have concluded that $R_{d,d',\nu}$ is not a constant multiple of a square polynomial modulo $p$ when $p$ is a prime such that $p|q$, $p\nmid dd'$ and $p > \max_{1 \leq i \leq k} h_i$, Lemma \ref{Weil} yields
\begin{align}\label{apwe}
\sum_{m \leq \frac{x}{dd'a_1[a_2^*,\dots,a_k^*]}} \chi(R_{d,d',\nu}(m)) \ll(d,q)^{1/2}(d',q)^{1/2}\left(\frac{x}{dd'a_1[a_2^*,\dots,a_k^*]q^{(1-\eps)/2}} +q^{(1+\eps)/2} \right).
\end{align}

Using Lemma \ref{lcm} and the condition $(a_i,a_j)|(h_i-h_j)$ from the outermost sum of (\ref{S1b}), we find that
\begin{align}\label{b1}
[a_2^*,\dots,a_k^*] \geq a_2\cdots a_k \prod_{1\leq i < j \leq k}(a_i,a_j)^{-1} \geq a_2\cdots a_k\left( \max_{1 \leq i \leq k} h_i \right)^{-\binom{k}{2}}.
\end{align}
We combine (\ref{apwe}) with (\ref{b1})
and since $d$ and $d'$ are coprime and the arithmetic function $\rho$ defined as $\rho(m)\vcentcolon=\{\nu\mod{m}:Q(\nu)\equiv0\mod{m}\}$ for all $m\in\N$ is multiplicative (by the Chinese Remainder Theorem), (\ref{Dp1}) gives
\begin{align}\label{forSigma'}
\begin{split}
\Sigma'\ll&\,\frac{x}{q^{(1-\eps)/2}a_1\cdots a_k}\bigg(\sum_{\substack{d\mid\prod_{p\leq z}p\\(d,a_1\cdots a_k)=1}}\frac{\rho(d)(d,q)^{1/2}}{d}\bigg)\bigg(\sum_{d'\mid a_1\cdots a_k\prod_{p\leq k}p}\rho(d')\mu(d')^2\bigg)\\
&+q^{(1+\eps)/2}\bigg(\sum_{d\leq z^u}\rho(d)\sqrt{d}\bigg)\bigg(\sum_{d'\mid a_1\cdots a_k\prod_{p\leq k}p}\rho(d')\sqrt{d'}\bigg).
\end{split}
\end{align}
Using the trivial inequalities $\tau(m),\rho(m)\leq m\, (m\in\N)$, we derive the bound 
$$\sum_{d\leq z^u}\rho(d)\sqrt{d}\ll z^{5u/2},$$
as well as the estimate
$$\sum_{d'\mid a_1\cdots a_k\prod_{p\leq k}p}\rho(d')\sqrt{d'}\ll (a_1\cdots a_k)^{3/2}\tau(a_1\cdots a_k)\leq (a_1\cdots a_k)^{5/2}.$$
Consequently, after recalling that $a_1,\ldots,a_k\leq z^u$ in the expression (\ref{S1b}) of $S$, (\ref{forSigma'}) implies that
\begin{align}\label{forSigma'2}
\Sigma'\ll&\,\frac{x}{q^{1/4}a_1\cdots a_k}\bigg(\sum_{\substack{d\mid\prod_{p\leq z}p\\(d,a_1\cdots a_k)=1}}\frac{\rho(d)(d,q)^{\frac{1}{2}}}{d}\bigg)\bigg(\sum_{d'\mid a_1\cdots a_k\prod_{p\leq k}p}\rho(d')\mu(d')^2\bigg)+z^{\frac{5(k+1)u}{2}}q^{\frac{1+\eps}{2}}.
\end{align}

We previously noted that any two roots of the polynomial $R_{d,d',\nu}$ cannot be congruent modulo $p$. By an almost identical reasoning, for the primes $p > \max_{1 \leq i \leq k} h_i$ that do not divide the product $a_1\cdots a_k$, one can show that the $k$ roots $$\overline{a_1[a_2^*,\dots,a_k^*]a_i^{-1}}(a_1r^*+h_i^*)a_i^{-1}-1$$
of $Q$ modulo $p$ are distinct. Thus, $\rho(p)=k$ for every prime $p > \max_{1 \leq i \leq k} h_i$ which does not divide $a_1\cdots a_k$. Furthermore, if $p\mid a_1\cdots a_k$, then $Q$ reduces to a constant or first degree polynomial modulo $p$. Therefore, for the primes $p\mid a_1\cdots a_k$, we have that $\rho(p)\leq 1$. We are going to use these two observations to bound the sums at (\ref{forSigma'2}). Firstly, we have that
\begin{align}\label{forSigma'3}
\begin{split}
\sum_{\substack{d\mid\prod_{p\leq z}p\\(d,a_1\cdots a_k)=1}}\frac{\rho(d)(d,q)^{1/2}}{d}&\ll \prod_{\substack{\max_{1 \leq i \leq k} h_i<p\leq z\\p\,\nmid\, qa_1\cdots a_k}}\bigg(1+\frac{\rho(p)}{p}\bigg)\prod_{\substack{p\mid q,\,p\,\nmid\,a_1\cdots a_k\\p>\max_{1 \leq i \leq k} h_i}}\bigg(1+\frac{\rho(p)}{\sqrt{p}}\bigg) \\
&\leq\prod_{\substack{\max_{1 \leq i \leq k} h_i<p\leq z\\p\,\nmid\, qa_1\cdots a_k}}\bigg(1+\frac{k}{p}\bigg)\prod_{\substack{p\mid q,\,p\,\nmid\,a_1\cdots a_k\\p>\max_{1 \leq i \leq k} h_i}}(1+k) \\
&\ll \tau_{k+1}(q)(\log z)^k \ll q^{1/20}(\log z)^k.
\end{split}
\end{align}
We also have
\begin{align}\label{forSigma'4}
\sum_{d'\mid a_1\cdots a_k\prod_{p\leq k}p}\rho(d')\mu(d')^2&\ll \prod_{p\mid a_1\cdots a_k}(1+\rho(p))\leq \tau(a_1\cdots a_k)\leq \tau(a_1)\cdots\tau(a_k).
\end{align}

We now insert (\ref{forSigma'3}) and (\ref{forSigma'4}) into (\ref{forSigma'2}) to complete the estimation of $\Sigma'$. This way we specifically obtain
\begin{align}\label{b2}
\Sigma' \ll \frac{x(\log z)^k\tau(a_1)\cdots \tau(a_k)}{q^{1/5}a_1\cdots a_k} + z^{5(k+1)u/2}q^{(1+\eps)/2}.
\end{align}

\subsection{Evaluation of the ``error term" of $\Sigma$ in (\ref{D3})} Upon denoting the Big-Oh term of (\ref{D3}) by
\begin{equation}\label{Sig''}
\Sigma'' := z^u+\sum_{r \geq \frac{u}{2}}2^{-r}\sum_{\substack{\ell \leq \frac{x}{a_1[a_2^*,\dots,a_k^*]} \\ P^{-}(Q(\ell))>z_r \\ (Q(\ell),a_1\cdots a_k\prod_{p \leq k} p)=1}}\tau(Q(\ell))^{2},
\end{equation}
the next step of the proof is to bound $\Sigma''$. We basically achieve this with an application of Lemma \ref{H} to the divisor sum of $\Sigma''$.

The inequality $\tau(ab) \leq \tau(a)\tau(b)\, (a,b \in \N)$ implies that
\begin{equation}\label{asubm}
\sum_{\substack{\ell \leq \frac{x}{a_1[a_2^*,\dots,a_k^*]} \\ P^{-}(Q(\ell))>z_r \\ (Q(\ell),a_1\cdots a_k\prod_{p \leq k} p)=1}}\tau(Q(\ell))^{2}\leq \sum_{\ell \leq \frac{x}{a_1[a_2^*,\dots,a_k^*]}}F(Q_1(\ell),\ldots,Q_k(\ell)),
\end{equation}
where
\begin{align*}
F(n_1,\dots,n_k) &\vcentcolon= \prod_{i=1}^{k}\tau(n_i)^2\mathds{1}_{P^{-}(n_i)>z_r}\mathds{1}_{(n_i,a_1\cdots a_k)=1},\quad \text{and}\\
Q_{i}(X) &\vcentcolon= \frac{a_1[a_2^*,\dots,a_k^*]}{a_i}X+\frac{a_1r^*+h_i^*-a_1[a_2^*,\dots,a_k^*]}{a_i},\quad \text{for }\quad i \in \{1,\dots,k\}.
\end{align*}

We now aim to apply Lemma \ref{H} to the right-hand side of (\ref{asubm}) as this will lead to a bound for $\Sigma''$. In order to use Lemma \ref{H}, we verify that all of its conditions are met. First of all, the polynomials $Q_1,\ldots,Q_k$ are irreducible and pairwise coprime because they are linear and have distinct roots. 
Now, suppose that $Q=Q_1\cdots Q_k$ has a fixed prime divisor $p$, or that we equivalently have that $\rho(p)=p$.
If $p|a_1\cdots a_k$, then the condition $(Q(\ell),a_1\cdots a_k\prod_{p \leq k} p)=1$ in (\ref{Sig''}) is not satisfied for any $\ell \in \N$, meaning that $\Sigma''=0$. If $p\nmid a_1\cdots a_k$, then $Q$ is a polynomial of degree $k$. This implies that $p = \rho(p) \leq k$, and so, in this case, $Q$ has a fixed prime divisor $p \leq k$. Then the condition $(Q(\ell),a_1\cdots a_k\prod_{p \leq k} p)=1$ in (\ref{Sig''}) is not met for any $\ell \in \N$ and $\Sigma''$ vanishes again. Hence, if $Q$ has a fixed prime divisor, then $\Sigma''=0$ and there is nothing to bound. Therefore, we may assume that $Q$ has no fixed prime factors. Finally, by combining $z^{30ku} = x^{\eps}$, which is an immediate consequence of the definition of $u$, with the inequalities $r^*\leq [a_2^*,\ldots,a_k^*]\leq a_2\cdots a_k$ and $a_1,\ldots,a_k\leq z^u$, we observe that
\begin{align*}
\frac{x}{a_1[a_2^*,\dots,a_k^*]} \gg\Vert Q \Vert^{\frac{1}{k}},
\end{align*}
where $\Vert Q \Vert$ denotes the sum of the absolute values of the coefficients of the polynomial $Q$.

Now that we finished verifying that Henriot's bound, namely Lemma \ref{H}, is applicable, we apply it to the right-hand side of (\ref{asubm}) and obtain that
\begin{align}\label{henr}
\begin{split}
\sum_{\substack{\ell \leq \frac{x}{a_1[a_2^*,\dots,a_k^*]} \\ P^{-}(Q(\ell))>z_r \\ (Q(\ell),a_1\cdots a_k\prod_{p \leq k} p)=1}}\tau(Q(\ell))^{2} \ll&\, \frac{x}{a_1[a_2^*,\dots,a_k^*]}\prod_{p\leq z}\left( 1-\frac{\rho(p)}{p} \right) \\
&\times \prod_{\substack{p|D \\ p \nmid a_1\cdots a_k}}\bigg( 1 + \sum_{\substack{\mu_1, \dots, \mu_k \in \{0,1\} \\ (\mu_1,\dots,\mu_k) \not= (0,\dots,0)}} \tau(p^{\mu_1})^2\cdots\tau(p^{\mu_k})^2 \bigg) \\
&\times \sum_{\substack{m_1,\dots,m_k \leq x \\ P^{-}(m_1\cdots m_k)>z_r \\ (m_1\cdots m_k,a_1\cdots a_k)=1}}\frac{\tau(m_1)^2\cdots\tau(m_k)^2}{m_1\cdots m_k},
\end{split}
\end{align}
where $D$ is the discriminant of $Q$. 

The value of $D$ can be easily evaluated by the roots of $Q$ and a simple calculation shows that
\begin{align*}
D = \frac{[a_2^*,\dots,a_k^*]^{k(k-1)}}{(a_1\cdots a_k)^{2(k-1)}}\prod_{1 \leq i < j \leq k}(h_i-h_j)^{2}.
\end{align*}
Since we are working with a fixed set of non-negative integers $h_1,\dots,h_k$, we note that
\begin{align}\label{henr2}
\prod_{\substack{p|D \\ p \nmid a_1\cdots a_k}}\bigg( 1 + \sum_{\substack{\mu_1, \dots, \mu_k \in \{0,1\}\\ (\mu_1,\dots,\mu_k) \not= (0,\dots,0)}} \tau(p^{\mu_1})^2\cdots\tau(p^{\mu_k})^2 \bigg) \ll 1.
\end{align}

We saw that
$\rho(p)=k$ for the primes $p > \max_{1 \leq i \leq k} h_i$ which do not divide $a_1\cdots a_k$. Therefore, by Bernoulli's inequality,
\begin{align}
\begin{split}
\prod_{p\leq z}\left( 1-\frac{\rho(p)}{p} \right) &\leq \prod_{\substack{\max_{1\leq i \leq k} h_i \leq p\leq z \\ p \nmid a_1\cdots a_k } }\left( 1-\frac{\rho(p)}{p} \right) \leq \prod_{\substack{\max_{1\leq i \leq k} h_i \leq p\leq z \\ p \nmid a_1\cdots a_k } }\left( 1-\frac{k}{p} \right) \\
&\leq \prod_{\substack{\max_{1\leq i \leq k} h_i \leq p\leq z \\ p \nmid a_1\cdots a_k } }\left( 1-\frac{1}{p} \right)^{k} \ll (\log{z})^{-k}\prod_{p | a_1\dots a_k}\left(1-\frac{1}{p}\right)^{-k} \\
&\leq (\log{z})^{-k}\prod_{j=1}^{k}\prod_{p | a_j}\left(1-\frac{1}{p}\right)^{-k}
= (\log{z})^{-k}\prod_{j=1}^{k}\left( \frac{a_j}{\varphi(a_j)} \right)^{k}. \label{henr3}
\end{split}
\end{align}

Furthermore,
\begin{align}
\begin{split}
&\sum_{\substack{m_1,\dots,m_k \leq x \\ P^{-}(m_1\cdots m_k)>z_r \\ (m_1\cdots m_k,a_1\cdots a_k)=1}}\frac{\tau(m_1)^2\cdots\tau(m_k)^2}{m_1\cdots m_k} \leq \Bigg( \sum_{\substack{m \leq x \\ P^{-}(m) > z_r}} \frac{\tau(m)^2}{m} \Bigg)^k \\
&\ll \prod_{z_r < p \leq x}\left( 1+\frac{4}{p}+O\left( \frac{1}{p^2} \right) \right)^k \ll \exp\left( \sum_{z_r < p \leq x} \frac{4k}{p} \right) \\
&\ll \left( \frac{\log{x}}{\log{z_r}} \right)^{4k}
= \left( \frac{25k}{25k-1} \right)^{4kr}\left( \frac{\log{x}}{\log{z}} \right)^{4k}.\label{henr4}
\end{split}
\end{align}

We insert the bounds (\ref{b1}), (\ref{henr2}), (\ref{henr3}), and (\ref{henr4}) into (\ref{henr}) and obtain
\begin{align*}
\sum_{\substack{\ell \leq \frac{x}{a_1[a_2^*,\dots,a_k^*]} \\ P^{-}(Q(\ell))>z_r \\ (Q(\ell),a_1\cdots a_k\prod_{p \leq k} p)=1}}\tau(Q(\ell))^{2} \ll \frac{x}{a_1\cdots a_k (\log{z})^k}\left( \frac{25k}{25k-1} \right)^{4kr}\left( \frac{\log{x}}{\log{z}} \right)^{4k}\prod_{j=1}^{k}\left( \frac{a_j}{\varphi(a_j)} \right)^{k}.
\end{align*}
But, we have that $25k-1>24k>8k/(\log 4-1)$, and so
\begin{align*}
\left( \frac{25k}{25k-1} \right)^{4kr} = \exp\left\{4kr\log\left( 1+\frac{1}{25k-1} \right) \right\} \leq \exp\left( \frac{4kr}{25k-1} \right)<2^re^{-r/2}.
\end{align*}
Hence,
\begin{align*}
\sum_{\substack{\ell \leq \frac{x}{a_1[a_2^*,\dots,a_k^*]} \\ P^{-}(Q(\ell))>z_r \\ (Q(\ell),a_1\cdots a_k\prod_{p \leq k} p)=1}}\tau(Q(\ell))^{2} \ll \frac{2^re^{-r/2}x}{a_1\cdots a_k (\log{z})^k}\left( \frac{\log{x}}{\log{z}} \right)^{4k}\prod_{j=1}^{k}\left( \frac{a_j}{\varphi(a_j)} \right)^{k},
\end{align*}
and then in virtue of (\ref{Sig''}), we obtain
\begin{align}
\begin{split}
\Sigma'' &\ll z^u+\frac{x}{a_1\cdots a_k (\log{z})^k}\left( \frac{\log{x}}{\log{z}} \right)^{4k}\prod_{j=1}^{k}\left( \frac{a_j}{\varphi(a_j)} \right)^{k}\sum_{r \geq u/2}e^{-r/2} \\
&\ll z^u+\frac{x}{a_1\cdots a_k (\log{z})^k}\left( \frac{\log{x}}{\log{z}} \right)^{4k}\prod_{j=1}^{k}\left( \frac{a_j}{\varphi(a_j)} \right)^{k}e^{-u/4}. \label{er}
\end{split}
\end{align}
\subsection{Completing the estimation of $S$ in (\ref{defS1})} We now combine (\ref{b2}) with (\ref{er}), and so (\ref{D3}) gives
\begin{align*}
\Sigma \ll&\, \frac{x(\log z)^k\tau(a_1)\cdots \tau(a_k)}{q^{1/5}a_1\cdots a_k} +z^{5(k+1)u/2}q^{(1+\eps)/2} \\
&+ \frac{xe^{-u/4}}{a_1\cdots a_k (\log{z})^k}\left( \frac{\log{x}}{\log{z}} \right)^{4k}\prod_{j=1}^{k}\left( \frac{a_j}{\varphi(a_j)} \right)^{k}.
\end{align*}
We insert this bound into (\ref{S1b}) and obtain
\begin{align}\label{endc1}
\begin{split}
S \ll&\, \frac{x(\log z)^k}{q^{1/5}}\sum_{\substack{a_1,\ldots,a_k\\P^+(a_1\cdots a_k)\leq z}}\frac{\tau(a_1)\cdots\tau(a_k)}{a_1\cdots a_k} +z^{(7k+5)u/2}q^{(1+\eps)/2}\\
&+\frac{xe^{-u/4}}{(\log{z})^k}\left( \frac{\log{x}}{\log{z}} \right)^{4k}\sum_{\substack{a_1,\dots,a_k \\P^{+}(a_1\cdots a_k)\leq z} } \frac{1}{a_1\cdots a_k }\prod_{j=1}^{k}\left( \frac{a_j}{\varphi(a_j)} \right)^{k}.
\end{split}
\end{align}
We now have that
\begin{align}\label{sm}
\begin{split}
&\sum_{\substack{a_1,\dots,a_k \\P^{+}(a_1\cdots a_k)\leq z} } \frac{1}{a_1\cdots a_k }\prod_{j=1}^{k}\left( \frac{a_j}{\varphi(a_j)} \right)^{k} \leq \Bigg( \sum_{\substack{ a \\P^{+}(a)\leq z} } \frac{1}{a}\left( \frac{a}{\varphi(a)} \right)^k \Bigg)^k \\
&\leq \prod_{p \leq z}\left( 1+\left( \frac{p}{p-1} \right)^k\frac{1}{p}+O\left( \frac{1}{p^2} \right) \right)^k \ll \exp\left( \sum_{p \leq z} \left( \frac{p}{p-1} \right)^k\frac{k}{p} \right) \\
&=\exp\left\{ \sum_{p \leq z} \left( \frac{k}{p}+O\left( \frac{1}{p^2} \right) \right) \right\} \ll (\log{z})^k.
\end{split}
\end{align}
In a similar fashion, one can show that
$$\sum_{\substack{a_1,\ldots,a_k\\P^+(a_1\cdots a_k)\leq z}}\frac{\tau(a_1)\cdots\tau(a_k)}{a_1\cdots a_k}\ll (\log z)^{2k}.$$
Upon recalling that $z=q^v$ and that $z^{30ku} = x^{\eps}$, the last bound, (\ref{endc1}), and (\ref{sm}) imply that
\begin{align}\label{part1}
S \ll \frac{x(\log z)^{3k}}{q^{1/5}} + xe^{-u/4}\left( \frac{\log{x}}{\log{z}} \right)^{4k}\ll \frac{xv^{3k}}{q^{1/6}} + xe^{-u/5}.
\end{align}
\subsection{The completion of the proof} Since $u=\eps V/(30kv)$, when we put the bound (\ref{part1}) together with the definition (\ref{defS1}) of $S$ and Lemma \ref{osz-sp}, we infer that
\begin{align}\label{bef4.1}
\sum_{n \leq x}\lambda_z(n+h_1)\cdots\lambda_z(n+h_k) \ll \frac{xv^{3k}}{q^{1/6}} + x\exp\left(-\frac{\eps V}{150kv}\right).
\end{align}

Now, if $V\geq 4\log\eta$, then from (\ref{vch}), it follows that $v=2$. Therefore, a combination of (\ref{bef4.1}) with (\ref{siegel}) and Lemma \ref{trans} yields
\begin{align}\label{case2i}
\begin{split}
\sum_{n \leq x}\lambda(n+h_1)\cdots\lambda(n+h_k) &\ll\frac{xV}{\eta}+\frac{x}{q^{1/6}}+x\exp\left(-\frac{\eps V}{300k}\right) \\
&\ll\frac{xV}{\eta}+x\exp\left(-\frac{\eps\sqrt{V\log\eta}}{150k}\right).
\end{split}
\end{align}

On the other hand, if $V<4\log\eta$, then according to the definition of $v$ in (\ref{vch}), we see that $v=\sqrt{V/(\log\eta)}<2$. So, by putting the estimate (\ref{bef4.1}) together with Lemma \ref{trans}, it follows that
\begin{align*}
\sum_{n \leq x}\lambda(n+h_1)\cdots\lambda(n+h_k) \ll&\,  \frac{x}{q^{1/6}}+x\exp\left(-\frac{\eps\sqrt{V\log\eta}}{150k}\right)\\
&+x\bigg(\frac{\log{\eta}}{V}\eta^{-\sqrt{\frac{V}{4\log{\eta}}}} +q^{-\sqrt{\frac{V}{\log\eta}}}\bigg)\bigg(V\sqrt{\frac{\log{\eta}}{V}}\bigg)^3 \\
\ll&\, \frac{x}{q^{1/6}}+x\exp\left(-c_0\sqrt{V\log\eta}\right)\\
&+x(V\log \eta)^{3/2}\exp\bigg(-\log{q}\sqrt{\frac{V}{\log{\eta}}} \bigg),
\end{align*}
where $c_0\in(0,\eps/(150k))$ is some constant.

According to (\ref{siegel}), we have that $xq^{-1/6}\ll x/\eta\ll xV/\eta$ and that $\log\eta\ll\log q$. Consequently, the last bound gives
\begin{align}\label{case2ii}
\sum_{n \leq x}\lambda(n+h_1)\cdots\lambda(n+h_k)\ll \frac{xV}{\eta}+x\exp\left( -c_1\sqrt{V\log{\eta}} \right),
\end{align}
where $c_1<c_0$ is a positive constant.

We now finish the proof of Theorem \ref{mainT} by merging the bounds (\ref{case2i}) and (\ref{case2ii}) from above.

\bibliographystyle{alpha}

\end{document}